\documentclass{amsart}
\usepackage{amsmath,amssymb,latexsym,amsthm}
\font\goth=eusm10

\newcommand\F{\mathcal F}
\newcommand\E{\mathcal E}
\newcommand\Ii{\hbox{\goth I}}

\newcommand\gH{\hbox{\goth H}}
\newcommand\Oc{\hbox{\goth O}}

\newcommand\Pj{\mathbb{P}}
\newcommand\ZZ{\mathbb{Z}}

\newtheorem*{rr}{Riemann-Roch formula (RR)}

\newtheorem{theorem}{Theorem}[section]
\newtheorem{proposition}[theorem]{Proposition}
\newtheorem{case}[theorem]{Case}
\newtheorem{definition}[theorem]{Definition}

\newtheorem{corollary}[theorem]{Corollary}

\newtheorem{remark}[theorem]{Remark}
\newtheorem{example}[theorem]{Example}
\newtheorem{lemma}[theorem]{Lemma}
\theoremstyle{plain}
\theoremstyle{definition}
\theoremstyle{remark}

\begin{document}
\large
\title{ACM bundles on a general quintic threefold.}
\author{Luca Chiantini, Carlo Madonna}
\email{chiantini@unisi.it, madonna@mat.uniroma2.it}
\subjclass{14F05}
\curraddr{Luca Chiantini: Dipartimento di Matematica, Via del Capitano, 15, 53100 SIENA, Italia} 
\curraddr{Carlo Madonna: Dottorato di Ricerca in Matematica,
Dipartimento di Matematica, Universit\`a degli Studi di Roma "Tor Vergata", Via della
Ricerca Scientifica, 00133 Roma, Italia}

\begin{abstract}
We give a partial positive answer to 
a conjecture of Tyurin (\cite {Tyu}). 
Indeed we prove that on a general quintic hypersurface of $\Pj^4$ 
every
arithmetically Cohen--Macaulay rank 2 vector bundle is infinitesimally rigid. 
\end{abstract}
\maketitle

\section{Introduction} \label{intro}
In this paper we study indecomposable vector bundles {\it without intermediate cohomology} on a smooth projective threefold $X \subset \Pj^n$ with Picard group generated over $\ZZ$ by an hyperplane section $H$. 

Let us recall the following:

\begin{definition}
Let $\E$ be an indecomposable rank $k$ vector bundle on  a smooth projective threefold $X$ as above. We say that $\E$ is an {\it arithmetically Cohen-Macaulay (ACM) bundle} if $h^i(\E(nH))=0$ for  $i=1,2$ and for any $n \in \ZZ$.
\end{definition}

In a previous paper (\cite{Mad1}) the second author showed a relation between the invariants (up to twist) of indecomposable rank $2$ ACM bundles on hypersurfaces of $\Pj^4$ and other threefolds of index one.\smallskip

Of course, the existence of ACM bundles is linked with the existence of some arithmetically Cohen--Macaulay curves, via Serre's celebrated correspondence between rank $2$ bundles on threefolds and subcanonical cuves. 

Remind that a projective, locally Cohen--Macaulay variety $X$ is {\it subcanonical} when the dualizing sheaf $\omega_X$ is $\Oc_X(eH)$ for some integer $e$. 

If the threefold $X$ is subcanonical itself (as complete intersection  are), then subcanonical curves $C$ on $X$ arise as $0$-loci of global sections of rank $2$ bundles $\E$ on $X$ and there is the natural exact sequence:
$$0\to \Oc \to \E\to \Ii_C(c_1(\E))\to 0\eqno {(0)}$$
where $\Ii_C$ is the ideal sheaf of the curve on $X$. Furthermore $C$ is ACM exactly when the bundle $\E$ is. $\E$ is decomposable if and only if $C$ is the intersection of $X$ with two hypersurfaces of $\Pj^n$.   

It follows that the existence and the behaviour under deformation of indecomposable rank $2$ bundles on $X$ is strictly linked with the problem of describing curves $C\subset X$, which are not complete intersection {\it on $X$}.\smallskip

Smooth Calabi--Yau threefolds $X$ which, in our terminology, are subcanonical with $\omega_X=\Oc_X$, received recently an increasing interest because of their connection with other fields of Mathematics. In particular, the study of curves on such threefolds (apart from the obvious complete intersection ones) and of their deformations  and, consequently, the study of related rank two vector bundles, was recently considered in the literature (see e.g. \cite {J-K} or \cite {Tyu}). 

The first examples of Calabi--Yau threefolds with Picard group generated by the hyperplane section, are smooth quintic threefolds in $\Pj^4$. In \cite{Tyu}, Tyurin conjectured that all rank 2 stable bundles $\E$ on a {\it general} quintic threefolds are {\it infinitesimally rigid}, that is the cohomology module $H^1(\E\otimes\E^{\vee})$, which represents the local deformation functor, vanishes. The conjecture of course imply that the Moduli spaces of these bundles on $X$ is a discrete set of points.
\smallskip

On the other hand, in a private communication, G. Ottaviani pointed out to us the following:

\begin{example} \cite{Ott2} \rm
Call $F$ the indecomposable Horrocks-Mumford rank $2$ bundle on $\Pj^4$ and call $\E$ the restriction of $F$ to a general quintic threefold $X$.\par
$F$ is not rigid in $\Pj^4$ (see \cite{H-M} for the properties of $F$). One has $h^0(F)>0$, $h^0(F(-1))=0$ and the cohomology of $F$ is shown in the table at p. 74 of \cite{H-M}.\par
Using the identification of Pic($X$) with $\ZZ$ and the exact sequence:
$$ 0\to F(-5) \to F \to E\to 0 \eqno {(1)}$$ 
one sees that $c_1(E)=5$ and $h^0(E(-1))=h^1(F(-6))=0$, hence $E$ is stable.\par
Furthermore, one computes that $H^1(F\otimes F^{\vee}) = 24$ and $$H^1(F\otimes F^{\vee}(-5)) = 0$$
(see e.g. \cite{D-S} p.218), hence by $ 0\to F\otimes F^{\vee}(-5)\to F\otimes F^{\vee}\to E\otimes E^{\vee}\to 0$
one obtains $H^1(E\otimes E^{\vee})>0$.\par
Thus $E$ is a counterexample to Tyurin's conjecture (in fact it is not even actually rigid).
\end{example}

The non--rigidity of $E$ follows euristically from the remark that $F$ must change under the action of PGL(4), for every homogeneous rank 2 bundle splits. Since no elements of PGL(4) fix a general quintic threefold in $\Pj^4$, clearly also $E$ cannot be rigid.
Observe that a similar argument would work for every indecomposable rank 2 bundles over $\Pj^4$, restricted to $X$.\smallskip

We remark that the bundle $E$ of the previous example is not ACM. Indeed from sequence (1) and from the table of \cite{H-M} p.74, one computes $h^1(E(2))=5$, since $h^1(F(2))=5$ and $h^1(F(-3))=h^2(F(-3))=0$.   \smallskip

Tyurin's conjecture is still open for stable rank 2 bundles on a general quintic threefold, which are not restriction of bundles in $\Pj^4$. Among these bundles, there are ACM bundles, which we are going to study in this note. 

Our aim is to prove the following:

\begin{theorem} \label{thm:rigidity} All stable rank 2 ACM bundle on a general quintic threefold are infinitesimally rigid.
\end{theorem}

Our result is based on the classification of invariants of indecomposable ACM rank 2 bundles on a smooth threefold in $\Pj^4$, obtained by the second author  in \cite{Mad1}. It turns out that if the bundle $E$ is normalized so that $h^0(E)>h^0(E(-1))=0$, then only few possibilities are left for the Chern classes of indecomposable rank 2 ACM bundles. Furthermore one has a classification of curves arising as 0-loci of these bundles (see section 3).
Then we use the method introduced by Kleppe and Mir\'o--Roig in \cite{K-M} to understand the infinitesimal deformations of ACM subcanonical curves which arise in each case, and we extablish the rigidity directly.   \smallskip

Let us recall that a classification of ACM bundles and ACM subcanonical curves on Fano threefolds are studied in the literature, and in some situation their moduli spaces are described. This is the case, for example, of rank 2 ACM bundles on some Fano hypersurfaces of $\Pj^4$, which are related to a pfaffian description of forms (see \cite {Bea1}). We refer the reader to \cite{Ott} for the quadric threefold, to \cite{A-C} \cite{Bea2} \cite{M-T} and \cite{Dru} for the cubic threefold, to \cite{I-M} and \cite{Mad2} for the quartic threefold and more  generally to \cite {A-G} \cite{Mad3} \cite{Kno} \cite{Bea1} \cite{S-W} and \cite{B-G-S}. 

For bundles on Calabi--Yau threefolds, our general references are \cite{Tho} and \cite{Tyu}. \smallskip
Notice that our method gives a description of the invariants of all possible non complete intersection subcanonical ACM curves on a general quintic threefold. Unfortunately in some cases the problem of their existence is still open. \smallskip
 
\noindent\textbf{Acknowledgements.} 
The authors are glad for this opportunity to contribute to the celebration of the 60th birthday of Prof. Silvio Greco. In particular the first author is grateful to Silvio, who was one of the advisors of his thesis, teached him a lot of commutative algebra and algebraic geometry, and encouraged him at the beginning of his career. \bigskip

\section{Generalities} \label{S:gener} \par
We work in the projective space $\Pj^4$ over the complex field. We will denote with $\Oc$ the structure sheaf of $\Pj^4$.\smallskip

Let $X$ be a  general quintic hypersurface in $\Pj^4$. 

$X$ is smooth and we identify its Picard group with $\ZZ$, generated by a hyperplane section. We use this isomorphism to identify line bundles with integers. In particular, for any vector bundle $\E$ on $X$, we set $c_1(\E)\in \ZZ$ and we write $\E(n)$ for $\E \otimes \Oc_X(n)$. 

We have the following formulas for the Chern classes of twistings of $\E$:
$$ c_1(\E(n)) = c_1(\E)+2n $$ 
$$ c_2(\E(n))= c_2(\E)+5nc_1(\E)+5n^2.$$

Let us define the number:
$$b(\E)=b=\max \{ n \mid h^0(\E(-n)) \ne 0 \}.$$
We say that $\E$ is {\it normalized} when $b=0$. Of course, after replacing $\E$ with the twist $\E(-b)$, we may always assume that it is normalized. 

We say that a rank 2 vector bundle $\E$ on $X$ {\it splits} if it is
isomorphic to the direct sum of two line bundles.

We use the definition of stability given in \cite{O-S-S} p. 160. Since $\text{Pic}(X) \cong \ZZ$, in our notation we know that:

\begin{remark} \label{stab} A rank 2 vector bundle $\E$ is semi-stable if and only if $2b-c_1\leq 0$. It is stable if and only if the strict inequality holds.
\end{remark}

In particular, when $\E$ is normalized, then it turns out that $\E$ is (semi-) stable if and only if $c_1>0$ ($\geq$).

The number $2b-c_1$ is invariant by twisting i.e. for all $n\in\ZZ$:
$$2b-c_1=2b(\E(n))-c_1(\E(n)).$$
It measures the {\it level of stability} of $\E$.
 
If $b=0$, it is shown in \cite{Har} Remark 1.0.1 that $\E$ has some global section whose zero-locus $C$ has codimension $2$. 
$C$ is a {\it subcanonical} curve of degree $c_2(\E)$, whose canonical divisor is $\omega_C=\Oc_C(c_1(\E))$.\smallskip

Since $\omega_X$ is trivial, Serre's duality says that:
$$ h^3(\E(n)) = h^0(\E^\vee(-n))= h^0(\E(-c_1-n)).$$

Let us finally recall the following:

\begin{rr} Let $\E$ be a rank 2 vector bundle on $X$. Then
$$
\chi(\E(n))=\frac{5}{6}c_1^3+\frac52nc_1^2+\frac52 n^2c_1+\frac{10}6n^3 -\frac{c_1c_2}{2}-nc_2 +\frac{25}{6}c_1+\frac{25}3n. 
$$
\end{rr}
\bigskip

\section{ACM subcanonical curves on a quintic threefold}

 Let us recall the main result of \cite{Mad1}, rephrased in our situation:

\begin{theorem} \label{Mad}
Let be $\E$ a normalized rank 2 ACM bundle on $X$. If $\E$ is indecomposable, then: 
$$-3<c_1(\E)<5.$$ 
\end{theorem}

We present here a rough classification of curves arising as 0-loci of sections of indecomposable ACM bundles on a quintic threefold. These computations were announced in \cite{Mad3}. In the first three cases, the bundles are not stable, by remark \ref{stab}. We list them here for the sake of completeness.
\smallskip

{\it Let $\E$ be a stable ACM rank 2 bundle on a smooth quintic threefold $X$.}

\begin{case}\label{retta}
Assume that $c_1(\E)=-2$. Then $c_2(\E)=1$ and $\E$ has a section whose 0-locus is a line.
\end{case}
\begin{proof} By \cite{Har} Remark 1.0.1 we know that $\E$ has a section whose 0-locus $C$ is a curve. By assumtions we have $h^1(\E)=h^2(\E)=0$ and by Serre's duality $h^3(\E)=h^0(\E(2))$. The exact sequence (0) of the introduction here reads:
$$0\to \Oc_X \to \E\to \Ii_C(-2)\to 0$$
hence $h^3(\E)=h^0(\E(2))= 15$. Now just use (RR) to see that necessarily $c_2(\E)=1$, i.e. $C$ has degree 1, so it is a line.
\end{proof}

\begin{case}\label{conic}
Assume that $c_1(\E)=-1$. Then $c_2(\E)=2$ and $\E$ has a section whose 0-locus is a conic.
\end{case}
\begin{proof} As above we know that $\E$ has a section whose 0-locus $C$ is a curve, with the exact sequence:
$$0\to \Oc_X \to \E\to \Ii_C(-1)\to 0.$$
By assumtions we have $h^1(\E)=h^2(\E)=0$ and by Serre's duality $h^3(\E)=h^0(\E(1))= 5$. Then one computes $c_2(\E)=2$, i.e. $C$ has degree 2. \par
We do not know if $C$ is reduced or irreducible. On the other hand, from the exact sequence (0) we know  that $h^0(\Ii_C(1))=h^0(\E(2))-15$. Since $h^3(\E(2))=h^0(\E(-1))=0$ and $\E$ is ACM, one may compute $h^0(\E(2))$ using (RR). It turns out that $h^0(\Ii_C(1))=2$, hence $C$ is a plane curve of degree 2, i.e. a conic.
\end{proof}

\begin{case}\label{elliptic}
Assume that $c_1(\E)=0$. Then either:
\begin{enumerate}
\item $c_2(\E)=3$ and $\E$ has a section whose 0-locus is a plane cubic;
\item $c_2(\E)=4$ and $\E$ has a section whose 0-locus is a complete intersection space curve;
\item $c_2(\E)=5$ and $\E$ has a section whose 0-locus is a non--degenerate elliptic curve.
\end{enumerate}
\end{case}
\begin{proof} Here $\E$ is semi-stable, by remark \ref{stab}. As above we know that $\E$ has a section whose 0-locus $C$ is a curve, with the exact sequence:
$$0\to \Oc_X \to \E\to \Ii_C\to 0.$$
Again $h^1(\E)=h^2(\E)=0$ while $h^3(\E)=h^0(\E)$, so the Euler characteristic of $\E$ cannot determine $c_2(\E)$ in this case. On the other hand $h^0(\Ii_C(1))=h^0(\E(1))-5$ and one may use (RR) to compute $h^0(\E(1))$, since $h^3(\E(1))=h^0(\E(-1))=0$ and $\E$ is ACM. It turns out that $h^0(\Ii_C(1))=c_2(\E)-5$. Thus $\deg(C)=c_2(\E)\le 5$. Furthermore $\omega_C=\Oc_C$, so $C$ is not a line, thus it cannot be contained in more than 2 independent hyperplanes of $\Pj^4$; hence $\deg(C)=c_2(\E)\ge 3$. \par
If $c_2(\E)=3$, then $h^0(\Ii_C(1))=2$ and $C$ is a plane cubic.\par
If $c_2(\E)=4$, then $h^0(\Ii_C(1))=1$ and $C$ is a space curve. It is well known that any arithmetically Cohen--Macaulay subcanonical curve in $\Pj^3$ is complete intersection (see e.g. \cite{Har}). The invariants tell us then that $C$ is complete intersection of two quadrics in $\Pj^3$.\par
Finally when $c_2(\E)=5$, then $C$ is a non--degenerate, elliptic ACM quintic in $\Pj^4$.
\end{proof}

Let us turn our attention to stable bundles.

\begin{case}\label{canonica}
Assume that $c_1(\E)=1$. Then either:
\begin{enumerate}
\item $c_2(\E)=4$ and $\E$ has a section whose 0-locus is a plane quartic;
\item $c_2(\E)=6$ and $\E$ has a section whose 0-locus is a complete intersection space curve, of type (2,3);
\item $c_2(\E)=8$ and $\E$ has a section whose 0-locus is a non-degenerate (possibly singular) canonical curve of genus 5.
\end{enumerate}
\end{case}
\begin{proof} As above we know that $\E$ has a section whose 0-locus $C$ is a curve, with the exact sequence:
$$0\to \Oc_X \to \E\to \Ii_C (1)\to 0.$$
Here $h^3(\E)=h^0(\E(-1))= 0$ and one computes $c_2(\E)=8-2h^0(\Ii(1))$, so that $c_2(\E)$ is even and $4\le c_2(\E)\le 8$.\par
 If $C$ is degenerate, just as in the previous case it must be complete intersection, and one concludes exactly as above.\par
Assume that $c_2(\E)=8$, so that $C$ is non--degenerate; since $\omega_C=\Oc_C(1)$, then $C$ has arithmetic genus 5. \par
Observe that by (RR) and Serre's duality, $H^0(\Ii_C(2))=3$. If the $3$ quadrics are independent, then $C$ is complete intersection in $\Pj^4$.
\end{proof}

\begin{case}\label{bicanonica}
Assume that $c_1(\E)=2$. Then $c_2(\E)\leq 14$. Furthermore $\E$ has a section whose 0-locus $C$ is contained in a vector space of quadrics of dimension $14-c_2= 14 -\deg(C)$.
\end{case}
\begin{proof} This is the most difficult case, in which we can say few things about the curves associated to $\E$. We have here $h^3(\E)=0$ by assumptions and by Serre's duality, hence by (RR) $h^0(\E)=15-c_2$ and the first claim follows since we are assuming $h^0(\E)\geq 1$.\par
For the second claim, just use the exact sequence (0).
\end{proof}

Notice that as $\deg(C)$ decreases in the previous example, then $C$ must be contained in a huge linear system of quadrics. In principle, this gives a lower bound for $\deg(C)$, hence for $c_2(\E)$. The bound is easily found when $C$ is reduced and irreducible, but unfortunately in general we cannot assume this properties.\par

\begin{remark} \rm
Later (see proposition \ref{c1=2}) we shall see that in the previous case, necessarily $c_2(\E)\geq 11$.\end{remark}

\begin{case}\label{tricanonica}
Assume that $c_1(\E)=3$. Then  $c_2(\E)= 20$. Furthermore $\E(1)$ is generated by global sections, so it has a section whose 0-locus $C'$ is a smooth, irrreducible curve.
\end{case}
\begin{proof} As above we know that $\E$ has a section whose 0-locus $C$ is a curve, with the exact sequence:
$$0\to \Oc_X \to \E\to \Ii_C (3)\to 0.$$
Here $h^3(\E(-1))=h^0(\E(-2))= 0$, so the Euler characteristic of $\E(-1)$ vanishes. From (RR) one computes $c_2(\E(-1))=10$, so that $c_2(\E)=20$.\par
The bundle $\E(1)$ is regular in the sense of Castelnuovo--Mumford (see \cite{Mum}), for $h^3(\E(1)(-3))=h^0(\E(-1))=0$ and $\E$ is ACM. Thus one knows that $\E(1)$ has a section whose 0/locus $C'$ is a smooth curve. Since $C'$ is a smooth ACM curve, then it is also connected, hence irreducible.
\end{proof}

\begin{case}\label{pfaffiana}
Assume that $c_1(\E)=4$. Then $c_2(\E)=30$ and $\E$ has a section whose 0-locus is a smooth irreducible ACM curve of degree $30$, not contained in cubic hypersurfaces and whose ideal sheaf is generated by quartics.
\end{case}
\begin{proof}  $\E$ has a section whose 0-locus $C$ is a curve, with the exact sequence:
$$0\to \Oc_X \to \E\to \Ii_C (4)\to 0.$$
Since $h^3(\E(-1))=h^0(\E(-4))= 0$, from (RR) one computes $c_2(\E(-1))$ and finds that $c_2(\E)=30$.\par
The bundle $\E$ is regular in the sense of Castelnuovo--Mumford. Indeed $h^3(\E(-3))=h^0(\E(-1))=0$ and $\E$ is ACM. Then as in the previous case we conclude that $\E$ has a section whose 0-locus is smooth irreducible. The two final claims follows soon from the exact sequence above.
\end{proof}

We collect all the possible values of $c_1(\E)$ and $c_2(\E)$ in the following table:

$$\begin{matrix}
c_1(\E) & c_2(\E) & informations \cr 
- 2 & 1 & \text{line} \cr
-1 & 2 & \text{conic} \cr
 0 & \bigg\{ \begin{matrix} 3 \cr 4 \cr 5 \end{matrix} & \begin{matrix}
 \text{plane cubic} \cr
 \text{space curve c.i. type (2,2)} \cr
\text{elliptic non--degenerate} \end{matrix} \cr
1 & \bigg\{ \begin{matrix} 4 \cr 6 \cr 8 \end{matrix} & \begin{matrix}
 \text{plane quartic} \cr
 \text{space curve c.i. type (2,3)} \cr
\text{canonical non--degenerate} \end{matrix} \cr
2 & \leq 14 & \text{non--degenerate, }h^0\Ii_C(2)=14-c_2(\E) \cr
3 & 20 & \text{not contained in quadrics} \cr
4 & 30 & \text{smooth, irreducible, generated by quartics} \cr
\end{matrix}$$ 
\bigskip

\begin{remark}
The relative Chern classes uniquely determine the normalization
of $\E$, in the sense that a stable ACM bundle with $h^0(\E)>0$  is normalized if and only if its two Chern classes are in the list above.\par
\rm Indeed observe that in all the previous cases $h^0(\E)>0$ just by (RR). When $c_1(\E)= 1,2$ by stability $\E(-1)$ cannot have non-zero sections, hence $\E$ is normalized. When $c_1(\E)=3$ then $c_2(\E)=20$ and $c_2(\E(-1))=10$, which is impossible for a normalized ACM bundle with first Chern class 1.
Finally if $c_1(\E)=4$, $c_2(\E)=30$ and $\E$ is not normalized, then necessarily $\E(-1)$ has a section. One computes $c_1(\E(-1))=2$ and $c_2(\E)=15$ so $\E(-1)$ cannot be normalized, which is impossible since $\E$ is stable. 
\end{remark}
\bigskip

\section{Deformations of ACM curves and the rigidity theorem}

Curves arising as 0-loci of sections of ACM bundles on $X$ are subcanonical ACM curves in $\Pj^4$. In particular, such curves are {\it arithmetically Gorenstein} (\cite{B-E}). The resolution of the ideal sheaf of arithmetically Gorenstein curves in $\Pj^4$ is described by the following:

\begin{proposition}\label{Gor}
Let $C$ be an $e$-subcanonical, ACM curve in $\Pj^4$ and call $\Ii$ the ideal sheaf of $C$ in $\Pj^4$. Then one has a resolution:
$$0\to \Oc(-e-5)\to \oplus \Oc(-b_i)  \to \oplus \Oc(-a_i)\to \Ii\to 0 \eqno{(2)}$$ 
which is self-dual, up to twisting. Hence if one orders the $a_i$'s and the $b_i$'s so that $a_1\leq\dots\leq a_n$ and $b_n\leq\dots\leq b_1$, then:
$$ \forall i\quad a_i = b_i-e-5.$$
\end{proposition}
\begin{proof} see \cite{B-E}
\end{proof} 

Since a curve $C$ which is 0-locus of a section of an ACM bundle on $X$ is locally complete intersection, then the (embedded) normal bundle of $C$ in $\Pj^4$ is well defined. The cohomology of this bundle is computed from the following formula of Kleppe and Mir\'o--Roig: 

\begin{theorem}\label{KM}
Let $C$ be an $e$-subcanonical, ACM curve in $\Pj^4$. Let $N_C$ be its normal bundle in $\Pj^4$. Then, with the previous notation, one can compute $h^0(N_C)$ from the formula:
$$ h^0(N_C)= \sum_{i=1}^n h^0(\Oc_C(a_i))+\sum_{1\leq i\leq j\leq n} \binom{-a_i+b_j+4}4 -$$ 
$$-\sum_{1\leq i\leq j\leq n} \binom{a_i-b_j+4}4 - \sum_{i=1}^n \binom{a_i+4}4 \eqno{(3)}$$
\end{theorem}
\begin{proof} see \cite{K-M}
\end{proof} 

Let us settle the link between deformation of a rank 2 bundle $\E$ on the quintic threefold $X$ and deformations of the ACM subcanonical curve $C\subset X$ which arises as the 0-locus of a global section of $\E$. 

First, we observe that stable ACM bundles have the following property: 0-loci of different sections are different.

\begin{lemma}\label{uniq}
Let $\E$ be a normalized, stable,  ACM bundle of rank 2 on a smooth quintic threefold $X$
and let $s,s' \in H^0(\E)$ be two independent global sections of $\E$. Call $C, C'$ their 0-loci. Then (as schemes) $C\neq C'$.
\end{lemma}
\begin{proof}
Assume that $s'$ vanishes on $C$. Then $H^0(\E \otimes \Ii_C)$ has dimension at least 2, for it contains the independent sections $s$ and $s'$.
Tensoring sequence (0) with $\E^\vee$ one gets:  
$$ 0 \to \E^{\vee} \to \E \otimes \E^{\vee} \to \E \otimes \Ii_C \to 0 \eqno (4) $$
Since $\E$ is stable and normalized, then $c_1(\E)>0$ hence: 
$$h^0(\E^{\vee})=h^0(\E(-c_1(\E))=0$$ 
and $h^1(\E^{\vee})=0$, since $\E$ is ACM. It follows that $h^0(\E \otimes \E^{\vee}) \geq 2$, which is absurd since $\E$ is stable and hence simple.
\end{proof}

The previous result is in fact a special case of the following more general statement:

\begin{proposition}
Let $X$ be a smooth quintic hypersurface in $\Pj^4$ and let $\E$ and $\F$ be two normalized stable ACM bundles of rank 2. Let $s$ and $t$ be two global sections of $\E$ and $\F$, respectively. Call $C$ the 0-locus of $s$ and $C'$ the 0-locus of $t$. If $C=C'$ (as schemes), then there exist an isomorphism $\E \cong \F$ which carries $s$ to $t$. 
\end{proposition}
\begin{proof}
Clearly $\E$ and $\F$ have the same Chern classes since the schemes $C$ and $C'$ have the same numerical character.
Tensoring the exact sequence:
$$0 \to \Oc_X \to \E \to \Ii_C(c_1) \to 0$$
 by $\F^{\vee}$ one gets:
$$0 \to \F^{\vee} \to \E \otimes \F^{\vee} \to \Ii_C \otimes \F \to 0.$$
By assumptions, $t$ defines a section of $\F \otimes \Ii_C$, which is the image of an element in $H^0(\E \otimes \F^{\vee})$, since $h^1(\F^{\vee})=0$. It means that $t$  induces a morphism $\varphi:\F \to \E$. Replacing $\E$ by $\F$, it also follows the existence of a morphism $\psi:\E \to \F$ induced by $t$.
The two morphisms interchange $s$ and $t$, hence the composition is non zero.  Since both $\E$ and $\F$ are simple, both $\varphi$ and $\psi$ are invertible and we are done.
\end{proof}

Call $E_C$ the normal bundle of $C$ {\it on $X$} and $N_C$ the normal bundle of $C$ in $\Pj^4$. We have the exact sequence:
$$0\to E_C\to N_C\to\Oc_C(5)\to 0 \eqno (4)$$
furthermore Serre's correspondence (see e.g. \cite{O-S-S}) implies:
$$E_C = \E\otimes \Oc_C$$ 

\begin{proposition}\label{deform} $h^0(E_C)\geq h^0(\E)-1$ and the stable ACM bundle $\E$ is infinitesimally rigid if and only if $h^0(E_C)= h^0(\E)-1$. \end{proposition}
\begin{proof}
The first claim follows soon from the exact sequence:
$$0\to H^0(\Ii_C\otimes \E)\to H^0(\E)\to H^0(\Oc_C\otimes\E)=H^0(E_C)$$
since $h^0(\Ii_C\otimes \E)=1$ by lemma \ref{uniq}.\par\noindent
$\E$ is infinitesimally rigid when $h^1(\E\otimes \E^\vee)=0$. Return to sequence (4) above and take cohomology:
$$ H^1(\E^\vee) \to H^1(\E\otimes \E^\vee) \to H^1(\Ii_C(c)\otimes \E^\vee)\to H^2(\E^\vee) $$
where $c=c_1(\E)$. Since $\E$ is ACM, one gets  $H^1(\E\otimes \E^\vee) = H^1(\Ii_C(c)\otimes \E^\vee)$. \par\noindent
Tensoring the exact sequence $0\to \Ii_C\to \Oc_X\to \Oc_C\to 0$ with $\E$ and recalling that $\E(-c)=\E^\vee$, one gets in cohomology:
$$ 0\to H^0(\Ii_C(c)\otimes \E^\vee) \to H^0(\E) \to H^0(E_C) \to 
H^1(\Ii_C(c)\otimes \E^\vee)\to H^1(\E)=0$$
and by lemma \ref{uniq}, necessarily $H^0(\Ii_C(c)\otimes \E^\vee)$ is generated by $s$. Hence the formula in the statement is equivalent to $H^1(\Ii_C(c)\otimes \E^\vee)=0$.
\end{proof}

Next let us turn our attention to the Hilbert schemes in $\Pj^4$. \smallskip

Call $\gH'$ the Hilbert scheme of ACM curves in $\Pj^4$ with degree $c_2(\E)$ and genus $1+(c_1(\E)c_2(\E))/2$. 

\begin{proposition} All the points  of $\gH'$ which parametrizes curves $C$ arising as 0-loci of sections of the indecomposable ACM bundles determine a smooth open subset $\gH\subset\gH'$ such that all $Y\in\gH$ are ACM and satisfy $h^1(\Oc_Y(5))=0$.
\end{proposition}
\begin{proof} 
Arithmetically Gorenstein curves in $\Pj^4$ are unobstructed by \cite{Miro}. All these curves $C$ are $c$-subcanonical for some $c=c_1(\E)<5$. Hence they satisfies:
$$ h^1(\Oc_C(5))=h^0(\Oc_C(c-5))=0$$
since $c-5<0$ and $h^1(\Ii_C(c-5))=0$. We are done since the vanishing of $H^1(\Oc_Y(5))$ describes an open subset of $\gH'$, by semicontinuity, as the ACM condition does.
\end{proof}  

Let $\Pj = \Pj^{125}$ be the scheme which parametrizes quintic threefolds in $\Pj^4$.
In the product $\gH\times\Pj$ one has the incidence variety (i.e. the {\it Hilbert flag scheme}) 
$$I = \{(C,X): X\text{ is smooth and } C\subset X\},$$
with the two obious projections $p:I\to \gH$ and $q:I\to \Pj$.

\begin{corollary} \label{smooth}
$I$ is smooth and the map $q:I\to\Pj$ has smooth general fibers.
\end{corollary}
\begin{proof}
The fiber of $I$ over $Y\in\gH$ is $\Pj(H^0(\Ii_Y(5)))$. We know that $h^1(\Oc_Y(5))=0$, so $h^0(\Oc_Y(5))$ is just the constant computed by Riemann-Roch, moreover $H^1(\Ii_Y(5))=0$. It follows that $I$ is a projective bundle over $\gH$. \par
Since we work in characteristic 0, the second claim follows soon by the theorem of generic smoothness.
\end{proof} 

\begin{proposition} \label{fibraq} Assume that $\dim(I)\le 125+u$. Then either $q$ is not dominant or, for $(C,X)\in I$ general, the normal bundle $E_C$ of $C$ in $X$ satisfies $h^0(E_C)=u$.
\end{proposition}
\begin{proof} The previous corollary says that the fiber $I_X$ of $q:I\to\Pj$ over a general quintic threefold $X\in\Pj$ is smooth. When $q$ dominates, then $I_X$ is $u$-dimensional, furthermore it is smooth, by the previous corollary. Since $I_X$ is an open subset of the Hilbert scheme of curves in $X$ and it contains all 0-loci of ACM bundles, the claim follows. 
\end{proof}

\begin{remark}\rm The above quoted fact that $I_X$ corresponds to the Hilbert scheme of curves in $X$ follows soon from the definition of the Hilbert deformation functor. Infinitesimally, this is encoded in the sequence $0\to E_C\to N_C\to \Oc_C(5)\to 0$ (see \cite{Sern}).\par
In our setting  the situation is readily understood. Call $\Ii$ the ideal sheaf of $C$ in $\Pj^4$ and $\Ii_C$ the ideal sheaf of $C$ in $X$. The tangent space to $I$ at $(C,X)$ is equal to the product $H^0(E_C)\times H^0(\Ii(5))/(h)$, where $h$ is an equation of $X$. Notice that $H^0(\Ii(5))= H^0(\Ii_C(5))/(h)$. The tangent space to the fiber is the kernel of the map: 
$$H^0(N_C)\times H^0(\Ii(5))/(h) \to H^0(\Oc_X(5))$$
and it is equal to $H^0(E_C)$, as it follows from the diagram:
$$\begin{matrix}
 & & & &	0					 &     &   0           \cr
 & & & &	\downarrow				 &     & \downarrow    \cr
 & & & &	H^0(\Ii_C(5))			 &  =  & H^0(\Ii_C(5)) \cr
 & & & &	\downarrow				 &     & \downarrow    \cr
 & & & &	H^0(N_C)\times H^0(\Ii(5))/(h) & \to & H^0(\Oc_X(5)) \cr
 & & & &	\downarrow				 &     & \downarrow    \cr
0 & \to & H^0(E_C) & \to & H^0(N_C)		 & \to & H^0(\Oc_C(5)) \cr
 & & & &	\downarrow				 &     & \downarrow    \cr
 & & & &	0					 &     &   0           \cr
\end{matrix}$$
\end{remark}
\medskip

Recall that if $\E$ is a stable, normalized ACM bundle on the smooth quintic threefold $X$, then $h^3(\E)= h^0(\E^\vee)=0$, so $h^0(\E)$ can be computed by (RR):
$$ h^0(\E) =  \chi(\E)=\frac{5}{6}c_1^3-\frac{c_1c_2}{2}+\frac{25}{6}c_1. \eqno {(5)} $$

We collect together all the previous result and get our main formula:

\begin{theorem} \label{thm} Let $X$ be a general quintic threefold in $\Pj^4$ and let $\E$ be a stable, normalized ACM rank 2 bundle on $X$; write $c_1$ and $c_2$ for the Chern classes of $\E$. Call $C$ the 0-locus of a section of $\E$ and call $N_C$ the normal bundle of $C$ in $\Pj^4$. If:
$$ h^0(N_C)\leq  \frac{5}{6}c_1^3-c_1c_2+\frac{25}{6}c_1+5c_2\eqno {(6)}$$
then $\E$ is infinitesimally rigid. 
\end{theorem}
\begin{proof}
We always have $\dim(\gH)\leq h^0(N_C)$. The projection $p:I\to\gH$ has fibers at $C$ of dimension $h^0(\Ii_C(5))-1$, which is equal to $125-h^0(\Oc_C(5))$. Now by Riemann-Roch, $h^0(\Oc_C(5))=5c_2-(c_1c_2)/2$; indeed $C$ has degree $c_2$ and $\omega_C=\Oc_C(c_1)$, furthermore $c_1<5$, by theorem \ref{Mad} and 
$$h^1(\Oc_C(5))=h^0(\Oc_C(c_1-5))=h^1(\Ii_C(c_1-5))=0.$$
It follows:
$$\dim(I)\leq h^0(N_C)+125-5c_2+\frac{c_1c_2}2.$$
>From proposition \ref{fibraq} one gets that either the projection $q:I\to\Pj$ is not dominant or $h^0(E_C)\leq h^0(N_C)-5c_2+c_1c_2/2$. But the first case cannot hold, for we assumed that a general quintic threefold has a vector bundle as $\E$, i.e. a curve $C\in\gH$. Then if formula (6) holds, using the previous remark one gets:
$$h^0(\E)-1\geq h^0(N_C)-5c_2+\frac{c_1c_2}2\geq h^0(E_C).$$
Proposition \ref{deform} tells us that in fact the equalities must hold and $\E$ is infinitesimally rigid. 
\end{proof}

The proof of our main theorem follows from the list of all possible curves $C$ arising as 0-loci of sections of ACM bundles and computing $h^0(N_C)$, by means of 
theorem \ref{KM}.\medskip

{\bf Proof of theorem 1.3.}\par\smallskip 

\noindent{\bf Case $c_1(\E)=4$.} \par
We know that $c_2(\E)=30$ and $\E$ is generated by global sections. Let $C$ be a curve arising as 0-locus of a section of $\E$. Then we saw in the previous section that we may assume $C$ smooth; $C$ is not contained in cubics and the ideal sheaf $\Ii_C$ is generated by quartics. By \ref{Gor} we have an autodual resolution of the ideal sheaf $\Ii$ of $C$ in $\Pj^4$ of the type:
$$0\to \Oc(-9)\to \oplus \Oc(-b_i)  \to \oplus \Oc(-a_i)\to \Ii\to 0 $$ 
where necessarly $a_i=4$ for all $i$. Then by duality $b_i=5$ for all $i$. Since $C$ is 4-subcanonical, then $h^0(\Oc_C(4))=$ genus of $C= 61$. Since $C$ is ACM, it follows that $h^0(\Ii(4))=9$. The resolution is:
$$0\to \Oc(-9)\to \Oc(-5)^9  \to \Oc(-4)^9 \to \Ii\to 0. $$ 
Now using formula (3) one computes soon $h^0(N_C)=99 $ and the rigidity of $\E$ follows from formula (6).

\smallskip\noindent{\bf Case $c_1(\E)=3$.} \par
This is very similar to the previous one. We know that $c_2(\E)=20$ and $\E(1)$ is generated by global sections. Let $C$ be a curve arising as 0-locus of a section of $\E$. Then $C$ lies in no quadrics and the ideal sheaf $\Ii_C$ is generated by quartics; it follows that the minimal generators of the ideal of $C$ have degree 3 or 4. Since $C$ is 3-subcanonical, then $h^0(\Oc_C(3))=$ genus of $C= 31$. Since $C$ is ACM, it follows that $h^0(\Ii(3))=4$. Hence in the resolution of $\Ii$ we have $a_i=3$ for three values of $i$, so that by duality $b_i=4$ for three values of $i$. We do not know in principle how many minimal generators of degree 4 one has for $\Ii$, since it may depend on the syzygies among the cubics (and in fact it may vary). The resolution reads:  
$$0\to \Oc(-8)\to \Oc(-4)^a\oplus \Oc(-5)^4  \to\Oc(-3)^4 \oplus \Oc(-4)^a\to \Ii\to 0.$$ 
Now one just computes that the contributions of terms of degree 4 in formula (3) cancel, so whatever $a$ is, one computes $h^0(N_C)=74 $ and the rigidity of $\E$ follows from (6).

\smallskip\noindent{\bf Case $c_1(\E)=2$.} \par
This is as usual the most difficult case. Put $d=c_2(\E)$; we just know that $d\leq 14$. If  $C$ is a curve arising as 0-locus of a section of $\E$, then its ideal sheaf is generated by quartics, furthermore $C$ is non degenerate and $h^0(\Ii(2))=14-d$.
The resolution reads:
$$0\to \Oc(-7)\to \Oc(-3)^b\oplus\Oc(-4)^a\oplus \Oc(-5)^{14-d}\to $$
$$  \to\Oc(-2)^{14-d}\oplus \Oc(-3)^a\oplus\Oc(-4)^b\to \Ii\to 0.$$ 

\begin{proposition} \label{c1=2} $d\geq 11$\end{proposition}
\begin{proof} If $d<11$, then $C$ is contained is a family of quadrics of affine dimension at least 4. These quadrics cannot have a common hyperplane. Indeed otherwise the residue hyperplanes must be independent and meet at a point: since $C$ is locally complete intersection, the point cannot be a component of $C$ (embedded or not), so $C$ is degenerate, contradiction.\par\noindent
It follows that two general quadrics through $C$ meet in a surface and $C$ is contained in a complete intersection curve of type $(2,2,4)$. The resolution of the ideal sheaf $\Ii'$ of the residue curve $C'$ is computed with the mapping cone of the diagram:
$$\begin{matrix} 
0 &\to& \Oc(-8)&\to &\Oc(-4)b\oplus\Oc(-6)^2&  \to & \Oc(-2)^2\oplus \Oc(-4) \cr
 &    & \downarrow & & \downarrow & & \downarrow  \cr
0 &\to& \Oc(-7)&\to &\Oc(-3)^b\oplus\Oc(-4)^a\oplus \Oc(-5)^{x}&  \to & \Oc(-2)^{x}\oplus \Oc(-3)^a\oplus\Oc(-4)^b  \end{matrix}$$
(where $x=14-d$). One gets:
$$0\to \Oc(-6)^{12-d}\oplus\Oc(-3)^b\oplus \Oc(-4)^{a-1}\to \Oc(-3)^{14-d}\oplus\Oc(-4)^b\oplus \Oc(-5)^a \to$$
$$ \to\Oc(-1)\oplus \Oc(-2)^2\oplus\Oc(-4)\to \Ii'\to 0.$$ 
Then $C'$ is a degenerate curve of degree 5 and genus 3 (clearly not reduced irreducible). Since $C'$ lies in an ireducible quadric and its ideal is generated by quartics, we may link it again using a complete intersection of type $(1,2,4)$ and get a new curve $C''\subset\Pj^3$ of degree 3, whose ideal sheaf $\Ii''$ (in $\Pj^4$) has a resolution which can be computed by the mapping cone again. It begins with
$$  \Oc(-1)^{13-d}\oplus \Oc(-2)\oplus\Oc(-3)^{a-1}\oplus\Oc(-4)^{b+1}\to \Ii''\to 0$$ 
which is impossible when $d< 11$, for 3 independent hyperplanes determine a line.
\end{proof}

Let us go back to our computation of $h^0(N_C)$.\par
{\it If $d=c_2=11$} then the resolution of $\Ii$ reads:
$$0\to \Oc(-7)\to \Oc(-3)^b\oplus\Oc(-4)^a\oplus \Oc(-5)^3 \to $$
$$\to 
\Oc(-2)^3\oplus \Oc(-3)^a\oplus\Oc(-4)^b\to \Ii\to 0$$ 
and comparing the first Chern classes of the sheaves, one obtains $b=a+2$.
We do not know $b$ exactly: it depends on the intersections of the 3 quadrics. Nevertheless one can apply formula (3) and see that every term with $b$ cancels. The output is $H^0(N_C)=47$, which by (6) gives the infinitesimal rigidity of $\E$.\par\noindent
{\it If $d=c_2=12$} then $C$ is contained in two quadrics, hence in the resolution of $\Ii$ there is at most one syzygy of degree 3. If there are no syzygies of degree 3, then by duality the curve is generated by cubics and one computes that one cubic generator is enough. It turns out that $C$ is complete intersection of type $(2,2,3)$, so that $h^0(N_C) = 2h^0(\Oc_C(2))+h^0(\Oc_C(3))=50$ and the inequality (6) holds. If there is one syzygy of degree 3, then the resolution is:
$$0\to \Oc(-7)\to \Oc(-3)\oplus\Oc(-4)^2\oplus \Oc(-5)^2 \to\Oc(-2)^2\oplus \Oc(-3)^2\oplus\Oc(-4)\to \Ii\to 0$$ 
from which one computes by (3) again $h^0(N_C)=50$. In any case theorem \ref{thm} implies the infinitesimal rigidity of $\E$.\par
{\it If $d=c_2=13$} $C$ is contained in one quadric and no syzygies of degree 3 are allowed. The resolution is:
$$0\to \Oc(-7)\to \Oc(-4)^4\oplus \Oc(-5) \to\Oc(-2)\oplus \Oc(-3)^4\to \Ii\to 0$$ 
so one computes $h^0(N_C)=53$, the inequality (6) holds and $\E$ is infinitesimally rigid.\par
{\it If $d=c_2=14$} then there are no quadric generators. The resolution is:
$$0\to \Oc(-7)\to \Oc(-4)^7\to\Oc(-3)^7\to \Ii\to 0$$ 
so one computes $h^0(N_C)=56$, the inequality (6) holds and $\E$ is infinitesimally rigid.
\par\smallskip

\noindent{\bf case $c_1(\E)=1$.} \par
Put $d=c_2(\E)$. We know that $d=4,6,8$. Let $C$ be a curve arising as 0-locus of a section of $\E$.\par
If $d=4$ or $d=6$, then $C$ is a space curve, hence it is complete intersection; in the first case it is complete intersection of type $(1,1,4)$ and one computes 
$h^0(N_C)= 2h^0(\Oc_C(1))+h^0(\Oc_C(4))=20$, while in the second case $C$ is of type $(1,2,3)$ and one has  $h^0(N_C)= h^0(\Oc_C(1))+h^0(\Oc_C(2))+h^0(\Oc_C(3))=28$. In any event, the inequality (6) holds and $\E$ is infinitesimally rigid.\par 
If $d=8$ then $C$ is non degenerate. Since the resolution is minimal, syzygies may arise only in degree 3. By duality  $\Ii$ is generated in degree 3 and the resolution has the form:
$$0\to \Oc(-6)\to \Oc(-4)^a\oplus \Oc(-3)^b  \to\Oc(-3)^b \oplus \Oc(-2)^a\to \Ii\to 0.$$ 
One computes $a=h^0(\Ii(2))=h^0(\E(1))-1=3$ by (RR), since $h^3(\E(1))=0$. It is impossible to determine $b$ exactly. Even if $C$ is smooth, we have $b=0$ for general canonical curves (complete intersection of 3 quadrics) and $b>0$ for trigonal ones.
Nevertheless in the formula (3) of theorem \ref{KM} all the terms containing $b$ cancel, and one finally obtains the value:
$$h^0(N_C)=36=\frac{5}{6}c_1^3-c_1c_2+\frac{25}{6}c_1+5c_2$$
hence by formula (6) of the previous theorem, $\E$ is infinitesimally rigid.

All the cases are examined and the proof of our main theorem is concluded. \qed  

\begin{remark} \rm
Notice that in all the previous cases we get in fact an {\it equality} in formula (6). It turns out that $I$ has  dimension exactly 125 for all the quoted situations, so the map $I\to\Pj$ may be dominant and we might expect that a general quintic threefold $X$ has $c_1$-subcanonical curves of degree $c_2$ for all the values of $c_1$ and $c_2$ listed in the previous section.\par
This is true in many cases. When $C$ is complete intersection, then its existence in a general quintic threefold can be computed using the method of Kley (see \cite{Kley}). Also the existence of lines and conics is well known, and for elliptic curves we refer to \cite{Kley} again. In the case $c_1=4$ and $c_2=30$ the existence result follows by Beauville's paper \cite{Bea1}, where it shown that the existence of such bundles is equivalent to the pfaffian representation of the quintic and in which it proved that a general quintic threefold is pfaffian.\par
However there are still cases where the existence is not known: $c_1=2$ and $c_2=11,13,14$ (which seems to be the most difficult) and $c_1=3$, $c_2=31$. See also \cite{I-M}, \cite{Mad2} and \cite{Mad3} for discussions on this subject.
\end{remark}

\begin{remark} \rm Even when $\E$ is a non-stable ACM bundle, some rigidity statement holds. Namely in this case, if $\E$ is normalized then $c_1(\E)\leq 0$ so that $h^0(\E)=1$ and one may replace the rigidity of $\E$ with the rigidity of the curve arising as 0-locus of a section of $\E$. This rigidity holds for lines, conics and elliptic curves in a general quintic threefold (see \cite{Kley}).
\end{remark} 
   
\begin{remark} \rm Buchweitz, Greuel and Schreyer proved in \cite{B-G-S} that every smooth quintic threefold has some non rigid bundle without intermediate cohomology. Our main theorem proves that such bundles must have rank at least 3 (see also \cite{Tyu} and \cite{Mad3} for a wider discussion of non rigid bundles on Calabi-Yau threefolds).
\end{remark}

We believe that our methods could be used to understand the rigidity (and somehow the existence) of ACM rank 2 bundles on general hypersurfaces in $\Pj^4$ of degree $d\geq 6$. For these threefolds, which are of general type, the existence of ACM indecomposable bundles would be in fact unexpected (see \cite{Voi}).

\end{document}